\documentstyle[12pt]{article}
\newcommand{\e}{\label}
\newcommand{\r}[1]{(\ref{#1})}
\newcommand{\be}{\begin{equation}}
\newcommand{\ee}{\end{equation}}
\newcommand{\bn}{\begin{eqnarray}}
\newcommand{\en}{\end{eqnarray}}
\newcommand{\ba}{\begin{array}}
\newcommand{\ea}{\end{array}}
\newcommand{\ibb}[1]{\leavevmode\hbox{\kern.3em\vrule
     height 1.2ex depth -.3ex width .2pt\kern-.3em\rm#1}}
\newcommand{\C}{{\ibb C}}
\newcommand{\Ibb}[1]{ {\rm I\ifmmode\mkern
            -3.6mu\else\kern -.2em\fi#1}}
\newcommand{\R}{{\Ibb R}}
\newcommand{\Z}{Z\!\!\!Z}

\newcommand{\lra}{\longrightarrow}
\newcommand{\ra}{\rightarrow}

\newcommand{\tr}{\mbox{{\rm{Tr}}}}
\def\qintd#1#2{\int_{#1}^{#2}\kern-1.5em\raise2.2ex\hbox{${\textstyle
q}$}\ } 
\def\qint{\int\kern-0.9em\raise2.2ex\hbox{${\textstyle
q}$}\ }

\begin{document}
\begin{titlepage}
\begin{flushright}
\begin{tabular}{r}
Helsinki Institute of Physics preprint\\
{\bf HIP~1997-02/Th}\\
\ \\
hep-th/97002132
\end{tabular}
\end{flushright}

\begin{center}
\vspace*{1.0cm}

{\Large\bf Quasiclassical Limit in q-Deformed Systems, Noncommutativity 
and the q-Path Integral}

\vskip 1cm

{\large {\bf M. Chaichian}}

\vskip 0.2cm

High Energy Physics Division, Department of Physics,\\ 
University of Helsinki \\ 
{\it and}\\ 
Helsinki Institute of Physics, \\
P.O.Box 9 (Siltavuorenpenger 20 C),\\ 
FIN-00014 Helsinki, Finland

\vskip 0.5cm

{\large {\bf A.P.~Demichev}}

\vskip 0.2cm

Nuclear Physics Institute,\\
Moscow State University,\\
119899, Moscow, Russia

\vskip 0.5cm

{\large {\bf P.P.~Kulish}}

\vskip 0.2cm

St.Petersburg Branch of\\ 
Steklov Mathematical Institute,\\
Fontanka 27,\\
191011, St.Petersburg, Russia

\end{center}

\vspace{1 cm}

\begin{abstract}
\normalsize

Different analogs of quasiclassical limit for a q-oscillator which result
in different (commutative and non-commutative) algebras of ``classical'' 
observables are derived. In particular, this gives the q-deformed
Poisson brackets in terms of variables on the quantum planes. We consider the
Hamiltonian made of special combination of operators (the analog of even
operators in Grassmann algebra) and discuss q-path integrals constructed with
the help of contracted ``classical'' algebras.  

\end{abstract}
\end{titlepage}

\

\vspace{3cm}

\section{Introduction}

Quasi-classical limit (QCL) in quantum mechanics is very important both
from general and utilitarian points of view. First of all, QCL provides
connection of classical and quantum theories (see, e.g. \cite{LandauL})
and, in fact, the very definition of quantization process is essentially 
based on the notion of QCL \cite{Berezin}. On the other hand, in 
appropriate situations there is effective quasi-classical approximation 
for calculations in quantum mechanics.

As is well known, quantization can be considered as the (continuous)
deformation ${\cal O}_\hbar$ of an algebra ${\cal O}_0$ of observables (real
functions) on classical phase space, so that a deformed associative product
(called star-product or $\star$-product) becomes non-commutative
\cite{Moyal,BayenFFLS}. Then QCL
appears as a result of expansion of a $\star$-product in powers of the
deformation parameter $\hbar$ and in terms of the non-deformed commutative
product. In particular, commutators of elements of ${\cal O}_\hbar$ (with
respect to a $\star$-product) give in first approximation in $\hbar$ the
Poisson brackets for elements (functions) of ${\cal O}_0$.

A few last years more general class of algebra deformations 
\cite{Drinfeld,ReshetikhinTF,Woronowicz2} (see also, e.g. 
\cite{ChariP,ChaichianD-book} and refs. therein) has attracted great
attention. One of the simplest but very important and, in a sense, basic
example of deformed algebras is the so called q-oscillator algebra 
\cite{Macfarlane,Biedenharn,PuszW,ChaichianK}. In an appropriate basis
this algebra ${\cal B}_q$ is generated by the elements $b,b^+$ with the
commutation relations (CR)
\be
bb^+-q^2b^+b=1 \ ,\qquad q\in\R  \ .                       \label{1}
\ee
Here $q$ is the deformation parameter of a usual oscillator (in
this paper we will consider only real values of $q$). 
The rescaling of the operators $b\ra b/\eta,\ b^+\ra b^+/\eta,
\ \ \eta\in\R$
brings to the CR \r{1} the second free parameter \cite{KulishD} 
and this CR becomes
\be
bb^+-q^2b^+b=\eta^2 \ ,\qquad q\in\R \ .                       \label{2}
\ee
We shall denote the corresponding algebra as ${\cal B}_{q,\eta}$ (all 
${\cal B}_{q,\eta}$ for $\eta\neq 0$ are isomorphic to ${\cal B}_{q,1}
\equiv {\cal B}_{q}$).
The values 
$$
q^2=1\ ,\qquad \eta^2=\hbar\ ,                $$
give the algebra ${\cal B}_{1,\hbar}$ of the usual harmonic 
oscillator with the CR 
\be
aa^+-a^+a=\hbar\ ,                                     \e{3}
\ee
(for this specific values of the parameters we use the ordinary notation
for oscillator operators: $b,b^+ |_{q^2=1,\eta^2=\hbar}=a,a^+$). 

The value $\hbar =0$ gives CR for the coordinates on classical
commutative phase space. Transition from \r{3} to the commutative case
corresponds to the usual quasi-classical limit in quantum mechanics.

We are going to consider two other limits for the deformation parameters
of ${\cal B}_{q,\eta}$; in other words, two other contractions of this
algebra: 
\begin{itemize}
\begin{enumerate}
\item 
The limit $\eta^2\ra 0,\ q=f(\eta):\ f(\eta)\ra 1$ as $\eta\ra 0$.

This leads to {\it commutative dynamics} in the classical limit. If, e.g. 
$q^2=1+\eta^2\gamma$, the CR (\ref{2}) becomes (cf. \cite{KulishD})
\be
[b,b^+]=\eta^2(1+\gamma b^+ b)\ ,                \e{3a}          \ee
and we expect to obtain in this limit classical dynamics with commutative
observables, curved phase
space and the symplectic form 
\be
\omega (\bar z,z)= i(1+\gamma \bar z
z)^{-1}dz\wedge d\bar z\ ,                   \e{3b}
\ee
($\bar z,z$ are coordinates on the classical phase space). In this case it
is worth to identify the parameter $\eta^2$ with the Planck constant
$\hbar$. 
\item 
The limit $\eta^2\ra 0$ with fixed  $q^2$. 
In this case one can expect to obtain q-classical dynamics with
{\it q-commuting variables} (coordinates on a q-plane)
\be
z\bar z=q^2\bar z z\ .                     \e{4}
\ee
\end{enumerate}
\end{itemize}

To study the first limit, we shall use the Bargmann-Fock representation for
the q-oscillator algebra and the corresponding symbol calculus in the space
with commuting variables but q-deformed differential and integral calculi.
Poisson bracket for this case have been introduced earlier
\cite{KulishD}, in the present paper we shall consider the series
expansion in the deformation constant $\eta^2$ of the $\star$-product.  Of
course, the resulting classical phase space is a usual manifold with ordinary
differential calculus (for $\gamma < 0$ the phase space proves to be
separated into two parts because of the degeneracy of the Poisson brackets
related to \r{3b}; this is the reason for the appearance of the non-Fock
representations for the q-oscillator with $q< 1$, see \cite{Kulish}).

The second limit is more unusual: the contracted algebra is also
non-commutative. 
The CR for the coordinates $\bar z,z$ on a q-plane can be considered as
the Weyl form of CR for a usual oscillator. Indeed, the operators $\bar
z,z$ can be represented in the form
\be
z=e^{\sqrt{\kappa}a}\ , \qquad\bar z=e^{\sqrt{\kappa}a^+}\ , \qquad
\kappa\in\R_+\ ,                                        \e{4a}
\ee
where $a,a^+$ satisfy \r{3}. As follows from the Baker-Hausdorff formula,
the operators $\bar z,z$ in this case obey \r{4} with
$q^2=e^{\kappa\hbar}$. Thus the contraction $\eta^2\ra 0$ (as well as the
contraction $q^2\ra 1$) of the q-deformed oscillator algebra leads to the
algebra of usual harmonic oscillator but with CR in Weyl form (on the
relation between Heisenberg and Weyl commutation relations, see
\cite{Faddeev}). 

The main aim of this paper is to study this limit (contraction to the Weyl
algebra) in details. Using
differential and integral calculi on Heisenberg-Weyl algebra (i.e. the
calculi on a q-plane) \cite{BaulieuF,KowalskiR,ChaichianDqpi} we shall
develop appropriate q-deformed symbol calculus and find the first
approximation (with respect to the parameter $\eta^2$) for the product in the
algebra of a q-oscillator in terms of Heisenberg-Weyl algebra. Therefore
we shall show that the symbol calculus can be applied to study 
a deformation (at least central one) of not only commutative algebras
(e.g., algebras of functions on classical phase spaces) but in the case of
deformations of non-commutative algebras as well. The expansions of the
product in the algebra ${\cal B}_{q,\eta}$ in powers of $\eta^2$ allows to 
define
q-deformed analog of a Poisson brackets and thus a kind of a q-deformed
``classical'' mechanics.

To define integral kernels of operators in terms of variables on a q-plane
one needs to consider arbitrary number of copies of ``coordinates'' on
q-planes. In other words, one must define extended algebra of arbitrary
number of pairs of operators $z$ and $\bar z$ with CR \r{4} with
differential and integral calculi on it. As a byproduct, with help of this
extended algebra we shall construct interesting Hamiltonian for
q-oscillator: with non-equidistant spectrum but harmonic dependence of the
operators $b(t),b^+(t)$ on time $t$. It is natural to call such a system
q-harmonic oscillator. This Hamiltonian belongs to the class of operators 
which resembles even combinations in Grassmann algebra.

The most important motivation for our consideration is necessity to
develop approximation methods for the study of systems with non-canonical
(deformed) commutation relations. One possible way to achieve this is
connected with q-deformed path integral
\cite{BaulieuF,ChaichianDqpi,ChaichianDqpi2}
heavily based on a symbol calculus. In the last Section we
shall discuss the meaning of q-deformed path integrals defined in terms of
non-commutative variables and its relation with q-deformed ``classical''
equation of motion.

\section{Bargmann-Fock representation for q-os\-ci\-llator algebra 
in terms of operators on quantum planes}

For the reader convenience we start from recalling of main formulas
concerning differential and integral calculi on a q-plane and
Bargmann-Fock (BF) representations for the q-oscillator
\cite{BaulieuF,KowalskiR,ChaichianDqpi,ChaichianDqpi2}. 

Two parametric differential calculus invariant with respect to quantum
Euclidean group $E(2)_{r,p}$ is defined by the following commutation
relations 
\begin{equation}
\begin{array}{ccc}
z\bar z = p\bar z z \ , & \qquad & \partial\bar\partial = p
\bar\partial\partial\ , \\
\qquad & \qquad & \qquad \\
dzd\bar z= - pd\bar z dz\ , &  \qquad &(dz)^2=(d\bar z)^2=0\ , \\
\qquad & \qquad & \qquad \\
zd\bar z = p(d\bar z) z\ , & \qquad & \bar z dz = p^{-1}(dz)\bar
z\ ,\\
\qquad & \qquad & \qquad \\
zdz = r^{-1}(dz) z\ , & \qquad & \bar z d\bar z = r(d\bar z)\bar
z\ ,\\
\qquad & \qquad & \qquad \\
\partial z=1+r^{-1}z\partial\ , & \qquad &
\bar\partial\bar z=1+r\bar z\bar\partial\ , \\
\qquad & \qquad & \qquad \\
\partial\bar z=p^{-1}\bar z\partial\ , &  \qquad &
\bar\partial z=p z\bar\partial\ , \\
\qquad & \qquad & \qquad \\
dz\bar\partial=p^{-1}\bar\partial dz\ , &  \qquad &
d\bar z\partial=p\partial d\bar z\ ,  \\
\qquad & \qquad & \qquad \\
d\bar z\bar\partial=r\bar\partial d\bar z\ , &  \qquad &
dz\partial=r^{-1}\partial dz\ .
\end{array}                              \label{2.2.20}
\end{equation}

In the standard BF representation the derivative $\bar\partial$ plays 
the role of
annihilation operator and a multiplication by the coordinate $\bar z$ 
plays the
role of creation operator. As is seen from (\ref{2.2.20}), q-deformed
$\bar\partial$ and $\bar z$ satisfy the CR (\ref{1}) for q-deformed
oscillator operators if one puts $r=q^2$. The second parameter $p$ can take
arbitrary value: this does not influence the CR (\ref{1}). On the other
hand properties of the quantum plane, which supposedly plays the role of a
classical phase space for the q-oscillator, strongly depends on the value of
parameter $p$: if $p=1$ the plane has commutative coordinates, if $p\neq
1$ the ``classical phase space'' proves to be non-commutative. 

Let us consider the construction of BF representation in two cases: 
\begin{itemize}
\begin{enumerate}
\item
non-commutative plane with $p=r=q^2$; 
\item
 commutative q-plane with $p=1,\ r=q^2$.
\end{enumerate}
\end{itemize}

Consider first the non-commutative space, $p=r=q^2$. 

To define a scalar product in BF
Hilbert space we need a q-analog of integration. In the case of usual
non-deformed BF representation the scalar product of two antiholomorphic
functions is defined with the help of exponential measure \cite{Bargmann}
\be
<g,f>=\int d\bar v dv\ e^{-\bar v v} \overline{g(\bar v)}f(\bar v) \ ,
                                               \label{2.2.20a}
\ee
($\bar v,v$ are commuting variables on ordinary complex plane).
In a q-deformed case corresponding integral is defined by the 
following postulates \cite{BaulieuF,KowalskiR,ChaichianDqpi}

\begin{itemize}
\begin{enumerate}
\item
Normalization
$$
I_{00}= \int d\bar zdz\ e_{q^2}^{-\bar z z} \ = 1\, ; $$
\item
Analog of the Stokes formula
$$
\int d\bar z dz\ \bar\partial f(\bar z,z)=\int d\bar z dz\ \partial
f(\bar z,z)=0\ ;$$
\item
Change of variables
$$
\int d\bar z dz\ e_{q^2}^{-a^2\bar z z}f(\bar z, z)=a^{-2}\int d\bar z dz\
e_{q^2}^{-\bar z z}f(a^{-1}\bar z,a^{-1}z)\ . 
  $$
\end{enumerate}
\end{itemize}
These postulates fix the form of the deformed exponent in the integrand 
and allow to prove the following result,
\begin{equation}
I_{mn}=\int d\bar z dz\ e_{q^2}^{-\bar z z}z^n\bar z^m = \delta_{mn}
[n;q^2]!\ ,\qquad n,m=0,1,2,...\, ,               \label{2.2.21}
\end{equation}
which is sufficient to compute the integral of any function which has
a power series expansion.

The BF representation is constructed in the Hilbert space $\cal H$
of anti-analytic functions $f(\bar z)$ on the q-plane with scalar
product of the form
\begin{equation}
<g,f>=\int d\bar z dz\ e_r^{-\bar z z} \overline{g(\bar z)}f(\bar z) \ ,
                                               \label{2.2.22}
\end{equation}
so that the monomials
\begin{equation}
\psi ^n (\bar z)=\frac{\bar z^n}{\sqrt{[n;q^2]!}}\ ,\qquad
[n;q^2]=\frac{q^{2n}-1}{q^2-1}\ ,                  \label{2.2.23}
\end{equation}
form the orthonormal complete set of vectors in $\cal H$. The creation
and annihilation operators are represented as coordinate and derivative
\begin{equation}
b^+=\bar z\ ,\qquad b=\bar\partial\ . \label{2.2.24}
\end{equation}
One can check that $b^+$ and $b$ are hermitian
conjugate of each
other with respect to the scalar product (\ref{2.2.22}).

Before going further we notice that a Hilbert space based on
operator-valued functions is rather unusual object in quantum mechanics.
But from general point of view the appearance of the scalar product 
\r{2.2.22} and calculations of operator averages with help of it does
not mean that we construct essentially new theory in compare with ordinary
quantum mechanics. To explain this, it is enough to remember that in general
case states of a system are represented in quantum mechanics by density
{\it operators} $\rho$, so that mean-value $\bar X$ of an observable
(operator) $X$ is defined by 
\be
\bar X = \tr (\rho X)\ .                                    \e{mean}
\ee
On the other hand, the scalar product 
\r{2.2.22} also can be understood as a trace of the integrand. Indeed, using
the representation of the algebra of operators on q-plane in 
$\ell^2$-space \cite{VaksmanK} spanned by vectors $|k\rangle_W,\ \ k\in\Z$:
\be
z|k\rangle_W=q^{-(k+1/2)}|k+1\rangle_W\ ,\qquad 
\bar z|k\rangle_W=q^{-(k-1/2)}|k-1\rangle_W\ ,           \e{z-rep}
\ee
one can check that for $q>1$ trace of the integrand in \r{2.2.21} does
exist and after the appropriate normalization coincides with the
value of the q-integral (right hand side in \r{2.2.21}) obtained by the pure
algebraic axiomatic way. Thus the calculation of operators mean-values with
help of the operator valued wave functions \r{2.2.23} again has the 
form \r{mean} with the integrand in \r{2.2.22} playing the role of a
density operator.

Now consider another possibility: plane with commuting
coordinates, $p=1,\ r=q^2$.
The BF representation of CR (\ref{1}) with commuting $z,\bar z$
variables and q-deformed differential and integral 
calculus is constructed in the
space of antiholomorphic functions of the form
$$
\psi_n=\frac{\bar z^n}{\sqrt{[n;q^2]!}}\ ,                   $$
which are orthonormalized with respect to the scalar product
$$
<g,f>=\int d\bar z dz\ e_{1/r}^{-q^2\bar z z}
\overline{g(\bar z)}f(\bar z) \ ,
$$
where the integral is defined by the relation 
\be
\int d\bar z dz\ e_{1/r}^{-q^2\bar z z}z^n\bar z^m = \delta_{mn}
[n;q^2]!\ .                                                \e{qint}
\ee
This integral is defined by the same postulate as in the non-commutative case.
Here we have used the so-called second q-deformed (basic) 
exponential function\index{q-exponent!second} \cite{Exton},
 $$
\exp_{1/q^2}\{x\}=\sum^\infty_{n=0}\frac{x^n}{[n;q^{-2}]!}=
\sum^\infty_{n=0}q^{n(n-1)}\frac{x^n}{[n;q^2]!}\ .
$$

This q-deformed integral can be interpreted in terms of Jackson integral, 
i.e. as sum over sites of a lattice. Thus in this
case there is quite ordinary quantum mechanical interpretation of the
scalar product and calculations of mean-values of observables.

\section{q-Oscillator as the deformation of commutative algebra}

As mentioned in Introduction, the adequate formalism to study the
quasi-classical limit is the symbol calculus 
\cite{Berezin,Moyal,BayenFFLS}. In the
case of oscillator-like systems, the most convenient symbol map is the {\it
normal} one, which put in correspondence to a regular function ${\cal
N}_A(\bar z,z)=\sum_{m,n}c_{mn}\bar z^mz^n$ the operator 
$A=\sum_{m,n}c_{mn}(b^+)^mb^n$, so that in each monomial of the operator,
creation operators are on the left of annihilation ones.
  
In the Bargmann-Fock representations this map can be made explicit 
using the operator kernels ${\cal A}$ of an operator
$A$ defined by the relation
$$
(Af)(\bar z_1) = \frac{1}{\eta^2}\int d\bar z_2 dz_2\ {\cal A}(\bar
z_1,z_2) {\cal E}(\bar z_2, z_2)f(\bar z_2)\ ,
$$
where the integration is understood in an appropriate sense (usual
integration for non-deformed oscillator or the integrations considered in
the preceding Section in q-deformed case) and ${\cal E}(\bar z,z)$ is the
function defined by Bargmann-Fock scalar products, e.g.
$$
\ba{rclcl}
{\cal E}(\bar z,z)&=&
e^{\scriptstyle{-\bar zz/\hbar}}\ ,&\quad&\mbox{in non-deformed case;}\\[3mm]
{\cal E}(\bar z,z)&=&
e_{\scriptstyle{1/q^2}}^{\scriptstyle{-q^2\bar zz/\eta^2}}\ ,
&\quad&\mbox{deformed BFR with}\ p=1,\ r=q^2\ ;\\[3mm]
{\cal E}(\bar z,z)&=&
e_{\scriptstyle{q^2}}^{\scriptstyle{-\bar zz/\eta^2}}\ ,
&\quad&\mbox{deformed BFR with}\ p=r=q^2\ .
\ea
$$

Of course, in cases of non-commutative variables one has to define CR
between different pairs and care about the order of factors in the
integrand. In non-deformed case the kernel of
an operator is expressed through the normal symbol by the relation
$$
{\cal A}(\bar z,z)=e^{\bar z z/\hbar}{\cal N}(\bar z, z)\ ,     $$
and the product of two operators corresponds to the convolution of their
kernels. This allows to find the explicit form of star-product of normal
symbols and its quasi-classical limit.

Now let us turn to the q-deformed cases.  
At first we will consider the algebra ${\cal B}_{q,\eta}$ 
in the limit:
$$\eta^2\ra 0,\qquad q^2(\eta)\ra 1$$
(the case 1 in Introduction). As we
expect to obtain classical dynamics with commutative observables, we 
express operators in terms of integral kernels with
commuting variables (but q-deformed differential and integral calculi with
$p=1,\ r=q^2$). Because of the same reason it is natural to consider the
parameter $\eta^2$ as the physical Planck constant and denote it as $\hbar$.

Manipulations analogous to those in the non-deformed case, 
give that the kernel
of product $A_1A_2$ of two operators $A_1,A_2$ with kernels ${\cal
A}_1,{\cal A}_2$, equals to the convolution
$$
{\cal A}_1\ast {\cal A}_2(\bar z,z)
=\frac{1}{\hbar}\ \int\ d\bar\xi d\xi e^{-q^2\bar\xi\xi/\hbar}_{1/r}
{\cal A}_1(\bar z,\xi){\cal A}_2(\bar\xi,z)\ ,
$$
constructed with help of the q-integral \r{qint}, 
and the relation between the normal symbol ${\cal N}$ and the
kernel ${\cal A}$ of an operator $A$ is
$$
{\cal A}(\bar z,z)=e^{\bar zz/\hbar}_r{\cal N}(\bar z,z)\ .   $$

So to understand how nontrivial (nonzero) $\hbar$ modifies the classical
(commutative)
multiplication we must calculate in the 
quasi-classical limit $\hbar\sim 0$ the expression
$$
{\cal N}_{{\cal A}_1}\star {\cal N}_{{\cal A}_2}=
 {\cal N}_{{\cal A}_1\ast{\cal A}_2}=
\frac{1}{\hbar}\int\  d\bar\xi d\xi\,{\cal N}_1(\bar z,\xi){\cal
N}_2(\bar\xi ,z)e_r^{\bar z\xi/\hbar}e_{1/r}^{-q^2\bar\xi\xi/\hbar}e_r^{\bar
\xi z/\hbar}
e_{1/r}^{-\bar zz/\hbar}\ .                   $$

In the usual classical case of the harmonic oscillator one would combine the
exponents, then shift the variables of integration and use, e.g., 
the steepest descent method to
evaluate the integral in $\hbar\ra 0$ limit. None of these methods can be
used in q-deformed case. Instead, we can use the 
q-exponent expansion and the q-integral property \r{2.2.21}.

Consider operators of the form $P_m=f_m(b^+)b^m$, where $m\in\Z_+$  and 
$f_m$ is an arbitrary polynomial. Any operator can be represented
as a sum (possibly infinite) of $P_m$. 
For the convolution of two such operators one has (after
trivial change of integration variables $\bar\xi,\xi\lra
\bar\xi/\sqrt{\hbar},\xi/\sqrt{\hbar}$)
$$
{\cal N}_{{\cal P}_m\ast{\cal P}_l}=\int  d\bar\xi d\xi\, f_m(\bar z)\xi^m
g_l(z)\bar\xi^l(\sqrt{\hbar})^{m+l}e_r^{\bar
z\xi/\sqrt{\hbar}}e_{1/r}^{-q^2\bar\xi\xi}e_r^{\bar \xi z/\sqrt{\hbar}}
e_{1/r}^{-\bar zz/\hbar}                         
$$
\be
={\sum_{s,n=0}^{\infty}}\int   d^2\xi f_m(\bar z)\xi^m
g_l(z)\bar\xi^l(\sqrt{\hbar})^{m+l-s-n}
\frac{\bar z^s\xi^s}{[s;q^2]!}\frac{\bar
\xi^nz^n}{[n;q^2]!}e_{1/r}^{-q^2\bar\xi\xi} e_{1/r}^{-\bar zz/\hbar}
\ .                                            \e{2.3.5}
\ee
For definiteness consider the case $l\leq
m$ (opposite case is dealt quite similarly).
Then from (\ref{2.2.21}) we have
$$
{\cal N}_{{\cal P}_m\ast{\cal P}_l}= f_m(\bar z) g_l(z) e_{1/r}^{-\bar
zz/\hbar}\sum_s \hbar^{l-s}\bar
z^sz^{m+s-l}\frac{[m+s;q^2]!}{[s;q^2]![s+m-l;q^2]!}  $$
\be
= f_m(\bar z) g_l(z) e_{1/r}^{-\bar zz/\hbar}\sum_s \hbar^{l-s}\bar
z^sz^{m+s-l}\frac{[s+m-(l-1);q^2]...[s+m;q^2]}{[s;q^2]!}\ .     \e{2.3.6}
\ee

To go further, let us note that the sum of the form
$$
\sum_{s=0}^\infty \bar z^sz^{s-l}\hbar^{l-s}\frac{1}{[s-k;q^2]!}\ ,       $$
where $k<l-1$, can be rewritten  up to $O(\hbar )$ terms as
\bn
\sum_{s=k}^\infty \bar z^sz^{s-l}\hbar^{l-s}\frac{1}{[s-k;q^2]!}&=&
\sum_{s=0}^\infty \bar
z^{s+k}z^{s+k-l}\hbar^{l-k}\hbar^{-s}\frac{1}{[s;q^2]!}\nonumber\\
&=&\bar z^k z^{k-l} \hbar^{l-k}e_r^{\bar zz/\hbar}\ .       \e{2.3.6a}
\en
Use of the summation theorem \cite{Exton} for q-exponentials with
commuting arguments 
\begin{equation}
e_q^Ae_{1/q}^B=\sum_{n=0}^\infty\frac{(A+B)^{(n)}_{q}}{[n;q]!}
\ ,\qquad AB=BA
\ ,                                      \label{2.3.gr2.3.39}
\end{equation}
where
$$
(A+B)^{(n)}_q=(A+B)(A+qB)\ldots (A+q^{n-1}B)\ ,                  $$
gives that this sum being inserted in (\ref{2.3.6}) instead of the sum there
would produce the expression of the order $\sim
O(\hbar^2)$. So in (\ref{2.3.6}) we must take into account the parts of the
last
factor of the form $\sim 1/[s-l;q^2]!$ and $\sim 1/[s-l+1;q^2]!$ only.

The identities
\be
[i+j;q^2]=q^j([i;q^2]-[-j;q^2])=q^j[i;q^2]+[j;q^2]\ ,           \e{2.3.8}
\ee
allow to present the ratio
\be
\frac{[s-l+1+m;q^2]\ldots [s+m;q^2]}{[s-l+1;q^2]\ldots [s;q^2]}\ , \e{2.3.7}
\ee
in the form
\bn
\lefteqn{q^{ml}+q^l[m;q^2]\frac{1}{[s-l+1;q^2]}
\sum_{i=0}^{l-1}q^{-i}}                          \e{2.3.11}\\
 &+& \mbox{terms with higher orders of}\ [s;q^2]\ \mbox{in denominators,}
\nonumber
\en
so that for \r{2.3.6} we have
\bn
\lefteqn{f_m(\bar z) g_l(z) e_{1/r}^{-\bar zz/\hbar}\sum_s \hbar^{l-s}\bar
z^sz^{m+s-l}\frac{[s+m-(l-1);q^2]...[s+m;q^2]}{[s;q^2]!}}\nonumber\\
&=& f_m(\bar z) g_l(z) e_{1/r}^{-\bar zz/\hbar}\sum_s \hbar^{l-s}\bar
z^sz^{m+s-l}\frac{q^{ml}}{[s-l;q^2]!} \nonumber\\
&+& f_m(\bar z) g_l(z) e_{1/r}^{-\bar
zz/\hbar}q^l[m;q^2]\nonumber\\
&\times&\left(\sum_{i=0}^{l-1}q^{-i}\right)\sum_s \hbar^{l-s}\bar
z^sz^{m+s-l}\frac{1}{[s-l+1;q^2]!}+O(\hbar^2)\ .               \e{2.3.12}
\en
The series over $s$ in this expression can be summed up to $e_r^{\bar
zz/\hbar}$, 
so that we finally obtain for the convolution, up to $\sim O(\hbar)$ terms
\be
{\cal N}_{{\cal P}_m\ast{\cal P}_l}=q^{ml}
{\cal N}_{{\cal P}_m}{\cal N}_{{\cal P}_l}
+\hbar  f_m(\bar z)z^{m-1} \bar z^{l-1}g_l(z)q^l[m;q^2]\sum_{i=0}^{l-1}q^{-i}
+O(\hbar^2)\ .         \e{2.3.13}
\ee

The factor $q^{ml}$ in the first term shows that in the limit $\hbar\ra 0$
the convolution does not reduce to usual commutative multiplication of
functions on``classical'' phase space. This means that there is no such a
classical limit for q-deformed oscillator with $q\neq 1$. From the other
hand, if $q=1+\gamma\hbar+O(\hbar^2)$, 
we have for (\ref{2.3.13}) in the same order in $\hbar$
$$
{\cal N}_{{\cal P}_m\ast{\cal P}_l}={\cal N}_{{\cal P}_m}{\cal N}_{{\cal
P}_l}(1+ml\gamma\hbar )+\hbar ml f_m(\bar z)z^{m-1} \bar z^{l-1}g_l(z)
+O(\hbar^2)   $$
\be
={\cal N}_{{\cal P}_m}{\cal N}_{{\cal P}_l}+\hbar  (1+\gamma\bar{z}z
)\partial{\cal N}_{{\cal P}_m}\bar\partial{\cal N}_{{\cal P}_l}
+O(\hbar^2) \ .                                              \e{2.3.14}
\ee
Here $\partial,\bar\partial$ denote, obviously, the usual (non-deformed)
derivatives.

The Poisson bracket derived from \r{2.3.14} has the form
\begin{eqnarray}
\{{\cal N}_1,{\cal N}_2\}_p&=&\lim_{\hbar\ra 0}\frac{i}{\hbar}
\left( {\cal N}_1\star{\cal N}_2-{\cal N}_2\star{\cal N}_1\right)\nonumber\\
&=&i(1+\gamma\bar z z)\left(\partial{\cal N}_1\bar\partial{\cal N}_2
-\partial{\cal N}_2\bar\partial{\cal N}_1\right)\ ,   \e{2.3.14a1}       
\end{eqnarray}
and, as we expected, it corresponds to the symplectic form \r{3b}.

\section{q-Oscillator as the deformation of Weyl algebra}

Now consider the second possibility for the quasi-classical limit mentioned 
in Introduction which leads to non-commutative 
``classical'' variables: 
$$\eta^2\ra 0\ ,\ \ \  q\ \mbox{is fixed.}$$
Here we use
the q-deformed BF representation with $p=r=q^2$ (see Section 2).
Again the action of any operator $A$ in the BF Hilbert space 
$\cal H$ based on non-commutative variables, can be
represented by its kernel ${\cal A}(\bar z,z)$
\begin{equation}
(Af)(\bar z_1) = \frac{1}{\eta^2}\int d\bar z_2 dz_2\ {\cal A}(\bar
z_1,z_2) e_r^{-\bar z_2 z_2/\eta^2}f(\bar z_2)\ ,
\label{2.3.gr2.3.26}
\end{equation}
where
\begin{equation}
{\cal A}(\bar z_1,z_2)=\sum_{m,n} A_{mn}
\frac{\bar z_1^m}{\sqrt{\eta^{2m}[m;q^2]!}}\frac{z_2^n}
{\sqrt{\eta^{2n}[n;q^2]!}}\ ,   \label{2.3.gr2.3.27}
\end{equation}
and we use here the q-integral \r{2.2.21}.
The special order of operator kernel and ``the integration measure''
(q-exponent) in \r{2.3.gr2.3.26} is convenient for the subsequent 
definition of convolution of
operators in the case of q-commuting ``classical'' variables. Another new
feature in this definition is that one more pair of q-commuting coordinates
is introduced. So we have to define the CR for coordinates on different
copies of q-planes. We postulate (cf. \cite{ChaichianDqpi}) that any
copies of
coordinates $\bar z_i,\ z_i \ (i=1,2,...)$ on q-planes have the
following CR:
\begin{eqnarray}
z_i\bar z_j=q^2\bar z_jz_i\ ,
\qquad \bar z_i\bar z_j=\bar z_j\bar z_i\ ,\nonumber \\
z_iz_j=z_jz_i\ ,                                          \label{2.3.gr2.3.28}
\end{eqnarray}
i.e. they do not depend on the indices which distinguish the copies.

Now we have to derive CR for differentials, derivatives and coordinates
for {\it different} pairs of variables. As usual we can do this using
the consistency requirement. For example, assume that CR for a derivative
$\partial_{z}$ and coordinate $\xi$ from another copy of variables is
a homogeneous one: $\partial_{z}\xi=a\xi\partial_z$, where $a$ is a constant
to be defined. Then, acting by both sides of this relation on a function
$f(z)$, we obtain $a=1$. Proceeding in this way, one comes to the following
CR for any two different pairs of variables $\{\bar z,z\}$ and
$\{\bar\xi,\xi\}$:

$$\ba{rclcrcl}
\partial_z\bar\xi & = & q^{-2}\bar\xi\partial_z\ ,&\qquad &\bar\partial_z\xi
&=& q^2\xi\bar\partial_z\ ,\\[2mm]
\partial_z\xi & = & \xi\partial_z\ ,&\qquad &\bar\partial_z\bar\xi &=&
\bar\xi\bar\partial_z\ ,                            \\[3mm]
dz\bar\xi & = & q^2\bar\xi dz\ ,&\qquad & d\bar z\xi 
&=& q^{-2}\xi d\bar z\ ,\\[2mm]
dz\xi & = & \xi dz\ ,&\qquad &d\bar z\bar\xi &=& \bar\xi d\bar z\ ,  \\[3mm]
dzd\bar\xi & = & q^2d\bar\xi dz\ ,&\qquad &d\bar zd\xi &=& q^{-2}d\xi 
d\bar z\ ,\\[2mm]
dzd\xi & = & d\xi dz\ ,&\qquad &d\bar zd\bar\xi &=& d\bar\xi d\bar z\ ,\\[3mm]
\partial_z\bar\partial_\xi & = & q^2\bar\partial_\xi\partial_z\ ,&\qquad
&\bar\partial_z\partial_\xi &=& q^{-2}\partial_\xi\bar\partial_z\ ,\\[2mm]
\partial_z\partial_\xi & = & \partial_\xi\partial_z\ ,&\qquad
&\bar\partial_z\bar\partial_\xi &=& \bar\partial_\xi\bar\partial_z\ ,\\[3mm]
\partial_zd\bar\xi & = & q^{-2}d\bar\xi\partial_z\ ,
&\qquad &\bar\partial_zd\xi
&=& q^2d\xi\bar\partial_z\ ,\\[2mm]
\partial_zd\xi & = & d\xi\partial_z\ ,&\qquad &\bar\partial_zd\bar\xi &=&
d\bar\xi\bar\partial_z\ .
\ea   $$

One can check that the definition of q-integral is consistent with these
CR for different variables if one requires that CR of the symbol $\int d\bar
z dz$ (map from q-plane to $\C$) are defined by (i.e. coincide with) the CR for
$d\bar z dz$.  This gives sense to the notation for the functional
(``definite integral'') on a q-plane.

Coefficients $A_{mn}$ in \r{2.3.gr2.3.27} can be expressed through 
the scalar product
\begin{equation}
A_{mn}=q^{2m(n+1)-2m}<\psi _m \mid A\mid\psi _n>\ .      \label{2.3.gr2.3.29}
\end{equation}

Consider an action of two operators $A_1$ and $A_2$ on arbitrary wave function
$f(\bar z)$
\begin{equation}
 A_2 A_1=\frac{1}{\eta^2}\int d\bar z_1 dz_1\ {\cal A}_2(\bar
z_2,z_1)  e_r^{-\bar z_1z_1/\eta^2}\frac{1}{\eta^2}\int 
d\bar z_0 dz_0 e_r^{-\bar
z_0 z_0/\eta^2}{\cal A}_1(\bar
z_1,z_0)f(\bar z_0)\ .                               \label{2.3.gr2.3.30}
\end{equation}
Using the CR for different pairs of q-variables, we obtain the formula for the
convolution of operator kernels
\be
{\cal A}_2\ast{\cal A}_1(\bar z,z)=\frac{1}{\eta^2}\int d\bar\xi d\xi \ {\cal
A}_2(q^{-2}\bar z,q^2\xi)  e_r^{-\bar \xi\xi/\eta^2} {\cal A}_1(\bar\xi,z) \ .
\label{2.3.gr2.3.30m}
\ee

Before going further, let us note that the evolution of an operator in 
quantum mechanics is defined by its commutator with the Hamiltonian of a 
system. The use of q-commutators which have the form
$$
[(b^+)^mb^n,(b^+)^kb^l]_q\equiv
\Bigl( (b^+)^mb^n\Bigr)\Bigl( (b^+)^kb^l\Bigr)
-q^{2(nk-ml)}\Bigl( (b^+)^kb^l\Bigr)\Bigl( (b^+)^mb^n\Bigr)\ ,
$$
(the q-factor being chosen accordingly to the homogeneous part of the
commutation relations \r{2}), contradicts
the Leibniz rule property of the time derivative in the Heisenberg
equation of motion, while the commutation
relations for general operators, made of $b^+,b$, are unnatural and
cumbersome. We recall that in somewhat analogous situation for fermionic
(anticommuting) operators, the problem is solved by the choice of
Hamiltonians from the even subalgebra of the complete fermionic algebra.
In our case, the operators which have natural commutation relations with
all other operators can be constructed with help of the number 
operator $N$ due to the equalities
\be
bq^{-2N}  = q^{-2}q^{-2N} b\ ,\qquad  b^+q^{-2N}  =
q^2 q^{-2N} b^+\, .                              \e{2.3.gr2.3.30m1}
\ee
It is easy to check that the monomials of the type
\be
q^{-2mN}(b^+)^mb^m\ ,\qquad m=1,2,...                 \e{add1}
\ee
have natural commutators with any other operators (i.e., the homogeneous
parts of the relations have no q-factors). Of course, the number
operator $N$ is not independent \cite{Kulish} and in the Fock
representation is connected with the creation and annihilation 
operators by the relation
\be
b^+b=\eta^2 [N;q^2];\qquad [N;q^2]\equiv\frac{q^{2N}-1}{q^2-1}\ ,\e{add2}
\ee
that is
$$
q^{2N}=1+\frac{1}{{\eta}^2}(q^2-1)b^+b\ .
$$
Thus, the operators \r{add1} read as
\be
\Bigl(1+\frac{1}{{\eta}^2}(q^2-1)b^+b\Bigr)^{-m}(b^+)^mb^m\ .     \e{add3}
\ee
However, there are a few reasons to treat the operator $q^{-2N}$ on a
footing distinguished from general operators made of $b^+,b$: {\it i)} the
relation \r{add2} is correct for the physically most important but still
particular Fock representation \cite{Kulish}; {\it ii)} the operators
\r{add3} are non-polynomial ones and, hence, finding of the normal form
for them is a quite complicated problem; {\it iii)} the most essential
reason in the context of our consideration is that the commutation
relations \r{2.3.gr2.3.30m1} are the {\it homogeneous} one and do not
depend on the algebra contraction parameter $\eta^2$ and on the 
representation. 

The latter fact inspires to introduce the notion of the q-deformed 
normal symbol. Let us denote for shortness
$$
K\equiv q^{-2N}\ ,
$$
and define a q-normal monomial
\be
N^p_{rs}=K^p(b^+)^rb^s\ ,             \e{2.3.gr2.3.30m3a}
\ee
with the corresponding q-deformed normal symbol
\be
{\cal N}^p_{rs}=K^p\bar z^rz^s\ .     \e{2.3.gr2.3.30m4}
\ee
The crucial difference from the usual normal symbol is that 
\r{2.3.gr2.3.30m4} is still an operator expression and that the operator 
$K\equiv q^{-2N}$ is the same as in the initial operator 
\r{2.3.gr2.3.30m3a}. We imply that the commutation relations of $K$ with 
$z,\bar z$ are the same as with $b,b^+$ 
\be
zK = q^{-2}K z\ ,\qquad
\bar zK = q^2 K \bar z\ .         \e{2.3.gr2.3.30m2}
\ee

>From \r{2.3.gr2.3.27} and \r{2.3.gr2.3.29} we can derive the expression for
the corresponding integral kernel
\be
{\cal A}^p_{rs}=q^{2(s(s+1)-r(p+1))}\bar z^r z^s
\exp_{1/q^2}\{q^{2(s-p+1)}\bar z z/\eta^2\}\ .
                                                    \label{2.3.gr2.3.35}
\ee

Now we are ready to calculate the normal symbol corresponding to the
product of two monomials in $K,b^+$ and $b$ in quasi-classical
approximation in analogy to the case of commuting variables. Consider the
convolution of two operator kernels
\bn
{\cal A}^p_{ab}\ast{\cal A}^t_{cd} & = &\frac{1}{\eta^2}\int
d\bar\xi d\xi q^{2(b(b+1)-a(p+1))}q^{-2a+2b}\bar
z^a\xi^b\exp_{1/q^2}\{q^{2(b-p+1)}\bar z\xi/\eta^2\} \nonumber\\[2mm]
&\times &\exp_{q^2}\{-\bar\xi\xi/\eta^2\} q^{2(d(d+1)-c(t+1))}\bar
\xi^cz^d   \nonumber\\
&\times &\exp_{1/q^2}\{q^{2(d-t+1)}\bar \xi z/\eta^2\}\ . \e{2.3m.1}
\en

Long but straightforward calculations of this convolution up to $O(\eta^2 )$
terms are analogous to those in the case of commuting variables. However, one
has to care about the order of all factors in the expressions. The result is
\bn
{\cal A}^p_{ab}\ast{\cal A}^t_{cd}&=&
q^{(b+d)(b+d+1)-(a+c)(p+t+1)}q^{2bc}q^{2t(a-b)}\bar
z^{a+c}z^{b+d} \nonumber\\[2mm]
&\times &\exp_{1/q^2}\{q^{2(b+d-p-t+1)}\bar zz/\eta^2\} \nonumber\\[2mm]
&+&q^{(b+d-1)(b+d)-(a+c-1)(p+t+1)}q^{2(b-1)(c-1)}q^{2t(a-b)}
[b;q^2][c;q^2]\nonumber\\[2mm]
&\times &\bar z^{a+c-1}z^{b+d-1}
\exp_{1/q^2}\{q^{2(b+d-p-t)}\bar zz/\eta^2\}+O(\eta^4) \ .         \e{2.3m.2}
\en

This convolution corresponds to the normal symbol
\bn
{\cal N}_{{\cal A}^p_{ab}\ast{\cal A}^t_{cd}}& = &
q^{2bc}q^{2t(a-b)}K^{p+t}\bar z^{a+c}z^{b+d} +\eta^2
q^{2(b-1)(c-1)}q^{2t(a-b)}[b;q^2][c;q^2] \nonumber\\[2mm]
& \times & K^{p+t}\bar z^{a+c-1}z^{b+d-1}+O(\eta^4) \ .
                    \e{2.3m.3}
\en

Introducing the star-product for normal symbols
\be
{\cal N}^p_{ab}\star{\cal N}^t_{cd}={\cal N}_{{\cal A}^p_{ab}\ast{\cal
A}^t_{cd}}\ ,                                                 \e{2.3m.3a}
\ee
we can present it in the form
\be
{\cal N}^p_{ab}\star{\cal N}^t_{cd}={\cal N}^p_{ab}{\cal N}^t_{cd}
+\eta^2({\cal N}^p_{ab}\partial_R)(\bar\partial{\cal N}^t_{cd})+O(\eta^4)
\ ,                                                            \e{2.3m.4}
\ee
where for convenience we have introduced right derivative
$\partial_R$  which
has the same CR as the left derivative $\partial$ but acts on functions from
the right. The same concerns its complex conjugate
$\bar\partial_R$. In particular,
$$
z^n\partial_R=[n;q^2]z^{n-1}\ ,\qquad 
\bar z^n\bar\partial_R=q^{-2(n-1)}[n;q^2]\bar
z^{n-1}\ .                                                  $$
Using \r{2.3m.4}, one can easily obtain the expression for q-commutator of
star-products of normal symbols in quasi-classical limit
\bn
\lefteqn{\left({\cal N}^p_{ab}\star{\cal
N}^t_{cd}-q^{2(bc-ad+p(d-c)+t(a-b))}{\cal N}^t_{cd}\star{\cal
N}^p_{ab}\right)}  \nonumber\\
&=&\eta^2{\cal
N}^p_{ab}\left(\partial_R\bar\partial-q^{2(b+c-1)}
\bar\partial_R\partial\right){\cal N}^t_{cd}+O(\eta^4)
\ .                                          \e{2.3m.5}
\en

An analogy with usual definition and procedure of physical quantization
inspires to define a q-deformed ``Poisson bracket'' as the factor in front
of $\eta^2/i$ in the rhs of \r{2.3m.5} in the limit $\eta^2\ra 0$
\be
\left\{{\cal N}^p_{ab},{\cal N}^t_{cd}\right\}_q=i{\cal
N}^p_{ab}\left(\partial_R\bar\partial-
q^{2(b+c-1)}\bar\partial_R\partial\right)
{\cal N}^t_{cd}\ .                                          \e{2.3m.6}
\ee
This means that in the limit
$\eta^2\ra 0$, the normal symbol of q-commutator
\be
[N^p_{ab},N^t_{cd}]_q=N^p_{ab}N^t_{cd}-
q^{2(bc-ad+p(d-c)+t(a-b))}N^t_{cd}N^p_{ab}  \ ,                \e{2.3m.7}
\ee
of two monomials of the form \r{2.3.gr2.3.30m3a} divided by $\eta^2/i$ is equal
to the q-Poisson bracket
\be
\lim_{\eta^2\ra 0}\frac{i}{\eta^2}{\cal N}_{[.,.]_q}=\{.,.\}_q
\ .                                  \e{2.3m.7a1}      
\ee

There is a special class of monomials, namely $N^a_{aa},\ \forall
a\in\Z$, for which the 
q-commutator \r{2.3m.7} becomes a usual commutator for any
second monomials entering the q-commutator
$$
[N^a_{aa},N^t_{cd}]_q=N^a_{aa}N^t_{cd}-N^t_{cd}N^a_{aa}\ ,
\quad \forall a,b,d,t\ .                                $$
As we discussed already, if we want to consider the 
Heisenberg equations of motion for q-oscillator
operators with usual time variables and, hence, with time derivative
satisfying  usual Leibniz rule, we have to use the
 Hamiltonians (time shift
generators)  of the type $N^a_{aa}$, since only the 
commutator acts on products of operators according to the Leibniz rule. 
Of course, this is true as we would like to
consider the classical limit of the type described 
above, so that the 
q-commutator corresponds to the q-Poisson brackets \r{2.3m.6}. 
If we would not care about
any limit at all, we could 
consider usual Heisenberg equations (with a commutator in its right hand
side) for any operators constructed from $b^+,b$.

Operators of the type $N^a_{aa}$ very much resemble even elements of a
superalgebra, in the sense that the latter have commutation (and not
anticommutation) relations with all other elements of the superalgebra. 

Thus the natural choice for Hamiltonian in our case is
\be
H_{qh}=N^1_{11}=\omega K b^+b\ ,\qquad \omega\in\R\ .  \e{2.3m.8}
\ee
In the q-oscillator basis $|n\rangle$ (see, e.g. 
\cite{ChaichianD-book}) 
this Hamiltonian is diagonal 
$$
H_{qh}\mid n\rangle =E_n\mid n\rangle \ ,        $$
where
\be
E_n=\frac{\omega}{q^2}[n;q^{-2}]\ .                  \e{2.3m.9}
\ee
The corresponding Heisenberg equations of motion 
prove to be very simple
\be
\frac{\partial b}{\partial t}=\frac{i}{\hbar}\left[H_{qh},b\right]
=-\frac{i\eta^2\omega}{\hbar}q^{-2}K b=
-\frac{i\eta^2\omega}{\hbar}bK\ ,                 \e{2.3m.10}
\ee
\be
\frac{\partial b^+}{\partial t}=\frac{i}{\hbar}\left[H_{qh},b^+\right]
=\frac{i\eta^2\omega}{\hbar}K b^+=
\frac{i\eta^2\omega}{\hbar}q^{-2}b^+K \ ,           \e{2.3m.10a}
\ee
with the obvious harmonic solution
\be
b(t)=e^{-iq^{-2}\eta^2\omega K t/\hbar}b(0)=
b(0)e^{-i\eta^2\omega K t/\hbar}\ ,                \e{2.3m.11}
\ee
\be
b^+(t)=e^{i\eta^2\omega K t/\hbar}b^+(0)=
b^+(0)e^{iq^{-2}\eta^2\omega K t/\hbar}\ .               \e{2.3m.11a}
\ee

``q-Quasiclassical'' equations of motion, i.e. the Heisenerg equations
\r{2.3m.10},\r{2.3m.10a} in the first approximation with respect to the
deformation parameter $\eta^2$, are defined by the q-deformed Poisson
brackets (cf. \r{2.3m.7a1})
\be
\frac{\partial z}{\partial t}=
\frac{\eta^2}{\hbar}\left\{H_{cl,qh},z\right\}_q
=-i\frac{\eta^2}{\hbar}q^{-2}\omega K z=
-i\frac{\eta^2}{\hbar}\omega zK\ ,                      \e{2.3m.13}
\ee
\be
\frac{\partial \bar z}{\partial t}=
\frac{\eta^2}{\hbar}\left\{H_{cl,qh},\bar z\right\}_q
=i\frac{\eta^2}{\hbar}\omega K z=
i\frac{\eta^2}{\hbar}q^{-2}\omega\bar zK \ .                    \e{2.3m.13a}
\ee
with the solutions of the form \r{2.3m.11},\r{2.3m.11a} 
as in the quantum case. Here we used
the q-normal symbol of $H_{qh}$ as classical Hamiltonian
\be
H_{cl,qh}={\cal N}_{H_{qh}}=\omega K\bar zz \ .        \e{2.3m.14}
\ee

If one puts $\eta^2=\hbar$ (but $q$ is still independent of $\eta^2=\hbar$; so
that $\kappa\sim\hbar^{-1}$ in \r{4a}), the equations \r{2.3m.13},\r{2.3m.13a}
prove to be true classical counterparts of \r{2.3m.10},\r{2.3m.10a}. In this
case there appears the interesting question about existence of classical
systems which have to be described by non-commutative algebras. In particular,
the possible area of applications might be the theory of stochastic systems
(with the parameter $q$ being related to physical parameters 
which are responsible for the stochastic behaviour of a system ). 

If the parameter $\eta^2$ is independent of $\hbar$ and, hence,
$\bar{z},\,z$ are considered as quantum mechanical operators, we observe the
amusing fact that the first order approximation to the q-oscillator with
respect to $\eta^2$ formally looks like q-generalization of classical
mechanics.

\section{Discussion}

We have developed deformed symbol calculus which allows to express
the product in the q-oscillator algebra as the power series in the
contraction parameter $\eta^2$ and in terms of the product of the contracted
algebra, i.e., in our case, commutative or Weyl algebra (algebra of 
variables on a quantum plane). The latter case is the non-commutative 
generalization 
of the usual quasi-classical approximation in quantum mechanics.

In the preceding section we already mentioned some analogy of our
consideration with the case of a superalgebra. Here we would like to add that 
the algebra describing spin-1/2 particles (in fact, algebra of Pauli matrices) 
also has different ``classical'' counterparts, either commutative or
anticommutative (Grassmann) algebras \cite{ChaichianD-gg}.

Interesting enough that though in general case contracted algebra is also
non-com\-mu\-ta\-tive, as well as initial one, equations of motion for its
operators are defined by a kind of Poisson brackets (q-deformed Poisson
brackets) and not by a commutator with a Hamiltonian. However, we have to
note that this q-Poisson bracket can not be represented in the coordinate
independent form contrary to the usual commutative case. Remind that in
the latter case Poisson brackets $\{\cdot,\cdot\}$ for functions 
$F$ and $H$ on a symplectic
manifold ${\cal M}$ with a symplectic 2-form 
$\Omega$ can be defined in the
coordinate-free form \cite{Arnold}
\be
\{F,H\}=\Omega (IdH,IdF)\ ,                        \e{pbr}
\ee
where the map $I:\ T^*{\cal M}_x\ra T{\cal M}_x$ 
(i.e., the map from 1-forms to vector
fields on ${\cal M}$) is defined by the relation
$$
\Omega(i,I\nu)=\nu(i)\ ,\qquad 
i\in T{\cal M}_x\ ,\nu\in T^*{\cal M}_x\ .
$$
One can consider the 2-form 
$\Omega_q=d\bar z\wedge dz$ as the symplectic
form on a q-plane but the peculiar form of the q-factor in q-Poisson
brackets \r{2.3m.6} prevents from a coordinate-free expression of the 
kind \r{pbr}. 

Another natural question which appears in connection with q-deformed calculus
is about possibility of construction and the meaning of q-deformed path
integral in terms of contracted but non-commutative algebras. The usual path
integral construction is heavily based on notions of a classical phase space.
This is one of reasons which emphasize the importance of the consideration in
Sections 3 and 4. Making use of the convolution of kernels for infinitesimal
evolution operators we have constructed such an integral for both cases of
the q-oscillator algebra contractions considered in the present paper
\cite{ChaichianDqpi,ChaichianDqpi2}. In particular, for the Hamiltonian
$H_{qh}$ (see \r{2.3m.8}) the evolution operator can be represented in the
form 
\be
U(\bar z,z;t^{\prime\prime}-t^\prime)= \int\left(\prod_t 
\frac{d\bar z(t)dz(t)}
{(1 + (1-q^2)\bar z(t) z(t))}\right)
 e_{1/r}^{q^2\bar z(t^{\prime\prime}) z(t^{\prime\prime})}
e^{\{iS_q\}}\ ,                                         \label{2.3m.17}
\ee 
where
\be
 S_q=\int^{t^{\prime\prime}}_{t^\prime}\left\{i \phi(\bar
z(t)z(t)) \bar z(t)\dot z(t)-\omega K\bar z(t)z(t)\right\} dt
\ ,                                                \label{2.3m.17a}
\ee 
$$
\phi(\bar z(t)z(t))=\sum^\infty_{r=0}\frac{q^{2r}}{(q^2(1-q^2))^{-1}+
q^{2r}\bar z(t)z(t)}\ ,                           $$
$\bar z(t),\ z(t)$ are different copies of coordinates on q-planes labeled
by the time parameter $t$.

This result shows that even for the Hamiltonian $H_{qh}$ in \r{2.3m.8} 
which has well
defined and simple dynamics in q-classical limit, the q-path integral
\r{2.3m.17}
(defined as the convolution of infinitesimal operator kernels) has a
complicated structure with nontrivial (non-flat) integration measure and
with the ``action'' (argument of the exponential in the integrand) which
does not lead to correct equations of motion \r{2.3m.13},\r{2.3m.13a}
through the formal application of usual principle of least action. 
The obvious reason for this is the
absence of the ordinary ``paths'' in a non-commutative space. 

However, we would like to note that this expression for q-deformed
path integral does not contradict the equations of motion 
\r{2.3m.13},\r{2.3m.13a}. Indeed, as we 
pointed out in Section 2, one can understand the ``q-integral'' over the
non-commutative variables $\bar z,z$ as the trace of an operator made of
only the combination $(\bar zz)$. 
Thus, in fact, the q-path integral \r{2.3m.17}, or rather its 
discrete
time approximation, can be understood as a sum over values of the operator
made of $\bar z(t_i)z(t_i)$ at different moments $t_i,\ i=0,\ldots,N$. The
distinction from a path integral over commutative phase space is that now
the operators $\bar z$ and $z$ do not have common eigenstates and, hence,
only one of them or their hermitian combination $\bar zz$ can have
definite eigenvalues. 

Thus in non-commutative case the path integral corresponds to summation
over ``reduced trajectories'' in the space of eigenvalues of maximal
set of commuting operators (in our case, just of one operator $\bar zz$).
This shows, in turn, that we may expect that evolution of only this operator
and its functions derived from the equations of motion 
\r{2.3m.13},\r{2.3m.13a} and from minimizing the
``action'' $S_q$ \r{2.3m.17a} (argument of the exponential in \r{2.3m.17}) 
must coincide. It is easy to see that this is indeed the case.

First of all, notice that from the equations of motion 
\r{2.3m.13},\r{2.3m.13a} it follows
that $\bar z(t)z(t)$ is an integral of motion (constant in time). On the
other hand, let us introduce the operators
$$
\rho(t)=\sqrt{\bar z(t)}\sqrt{z(t)}\ ,
\qquad \tau(t)=\sqrt{z(t)}\biggl(\sqrt{\bar z(t)}\biggr)^{-1}\ ,
$$
so that
$$
z=\tau\rho\ ,\qquad \bar z=\frac{1}{\sqrt{q}}\rho\tau^{-1}\ .
$$
Thus the ``action'' $S_q$ takes the form
\begin{eqnarray*}
S_q&=&\int_{t'}^{t''}\left(i\phi(\rho^2(t)/\sqrt{q})
\frac{1}{\sqrt{q}}\rho\tau^{-1}(\dot{\tau}\rho
+\tau\dot{\rho}) -\frac{\omega K}{\sqrt{q}}\rho^2(t)\right) dt \\
&=&\int_{t'}^{t''}\left(\frac{i}{\sqrt{q}}\phi(\rho^2(t)/\sqrt{q})
(\rho^2\dot{\nu}+\rho\dot{\rho}) -
\frac{\omega K}{\sqrt{q}}\rho^2(t)\right) dt \ ,\\
\end{eqnarray*}
where $\nu=\log\tau$. Making formal shift $\nu(t)\ra\nu(t)+\delta\nu(t)$
one obtains that the condition $\delta S_q=0$ implies
$$
\frac{d}{dt}\rho^2\phi\left(\rho^2(t)/\sqrt{q}\right)=0\ ,
$$
and hence $\rho(t)=const$.

Thus the evolution of the operator which defines ``trajectories'' in the
non-commutative case corresponds to the ``minimization'' of the q-deformed
action.

It is worthwhile to compare again this situation with the case of path integral
over Grassmann variables. There one also integrates, in fact, over even
combination of the variables (of the type $\bar zz$, as above). But additional
structures, analogous to the function $\phi(\bar zz)$, do not appear due to
nilpotence property of Grassmann variables. 

Generalization of our considerations to the case of multioscillators 
(with the final aim at a field theory construction) is straightforward;
similar treatments for other algebras and as well for the values of $q$ equal
to root of unity, remain an interesting problem for further investigations.

\vskip 10 mm

\centerline{{\bf Acknowledgments}}

A.P.D. thanks the Department of Physics and the Research 
Institute for High Energy
Physics, University of Helsinki, where part of this work has been done,
for hospitality; his work was partially supported by the INTAS-93-1630-EXT
and RFBR-96-02-16413-a grants. The work of P.P.K. was supported by the 
RFFI grant 96-01-00851.


\begin{thebibliography}{99}
\bibitem{LandauL}
L.D.Landau and E.M.Lifshitz, {Quantum Mechanics} (Pergamon Press, New
York, 1981)
\bibitem{Berezin}
F.A.Berezin, {\it Comm.Math.Phys.}, {\bf 40} (1975) 153.
\bibitem{Moyal}
J.E.Moyal, {\it Proc. Cambridge Phil. Soc.} {\bf 45} (1949) 99.
\bibitem{BayenFFLS}
F.Bayen, M.Flato, C.Fronsdal, A.Lichnerowicz and D.Sternheimer, 
{\it Ann. Phys.}
{\bf 111} (1978) 61; {\it Ann. Phys.} {\bf 111} (1978) 111.
\bibitem{Drinfeld}
V.G.Drinfel'd, in {\it Proc.of the International Congress of Mathematicians
(Berkley, 1986)} (American Mathematical Society, 1987), p.798.
\bibitem{ReshetikhinTF}
L.D.Faddeev, N.Yu.Reshetikhin and L.A.Takhtadjan, {\it Algebra i Analis}
{\bf 1} (1989) 178 (Transl.: {\it Leningrad Math. Jour.} {\bf 1} (1990)
193). 
\bibitem{Woronowicz2}
S.L.Woronowicz, {\it Comm. Math. Phys.} {\bf 111} (1987) 613.
\bibitem{ChariP}
V.Chari and A.Pressley, {\it A Guide to Quantum Groups}, (Cambridge Univ. 
Press, Cambridge, 1994) 
\bibitem{ChaichianD-book}
M.Chaichian and A.Demichev, {\it Introduction to Quantum Groups} (World
Scientific, Singapore, 1996).
\bibitem{Macfarlane}
 A. Macfarlane, {\it J. Phys.} {\bf A22} (1989) 4581.
\bibitem{Biedenharn}
L.C. Biedenharn, {\it J. Phys.} {\bf A22} (1989) L873.
\bibitem{PuszW}
W. Pusz and S.L. Woronowicz, {\it Rep. Math. Phys.}
             {\bf 27} (1989) 231.
\bibitem{ChaichianK}
M. Chaichian and P.P. Kulish, {\it Phys. Lett.} {\bf B234} (1990) 72.
\bibitem{KulishD}
P.P. Kulish and E. Damaskinsky, {\it J.Phys.} {\bf A23} (1990) L415.
\bibitem{Faddeev} 
L.D. Faddeev, {\it Lett.Math.Phys.} {\bf 34} (1995) 249.
\bibitem{BaulieuF}
L. Baulieu and E.G. Floratos, {\it Phys. Lett.} {\bf B258} (1991) 171.
\bibitem{KowalskiR}
K.Kowalski and J.Rembielinski, {\it J. Math. Phys.} {\bf 34} (1993) 2153.
\bibitem{ChaichianDqpi}
M.Chaichian and A.P.Demichev, {\it Phys. Lett.} {\bf B320} (1994) 273.
\bibitem{ChaichianDqpi2}
M. Chaichian and A.Demichev,
in  {\it Proceedings of the Clausthal Symposium on
``Nonlinear, Dissipative, Irreversible Quantum Systems''},  eds.  H.-D.
Doebner, V.K. Dobrev and P. Nattermann (World Scientific, Singapore,1995)
p.401.
\bibitem{Bargmann}
V.~Bargmann, {\it Proc. Nat. Acad. Sci. USA} {\bf 2} (1962) 199.
\bibitem{VaksmanK}
L.L.Vaksman and L.I.Korogodsky, {\it Dokl. Akad. Nauk SSSR} {\bf 304} (1989)
1036.
\bibitem{Exton}
H.Exton, {\it q-Hypergeometric functions and applications} (Ellis Horwood,
Chichester, 1983).
\bibitem{Kulish}
P.P. Kulish, {\it Theor. Math. Phys.} {\bf 86} (1991) 108.
\bibitem{ChaichianD-gg}
M.Chaichian and A.P.Demichev, {\it Phys. Lett.} {\bf A207} (1995) 23.
\bibitem{Arnold}
V.I. Arnold, {\it Mathematical Methods in Classical Mechanics} (Springer,
Berlin, 1978).

\end{thebibliography}
\end{document}